\documentclass[twoside,leqno,10pt]{amsart}
\usepackage{amsfonts}
\usepackage{amsmath}
\usepackage{amscd}
\usepackage{amssymb}
\usepackage{amsthm}
\usepackage{latexsym}
\usepackage{color}
\newtheorem{theorem}{Theorem}

\newtheorem{conjecture}{Conjecture}
\numberwithin{equation}{section}

\title{Variance of primes in short residue classes for function fields}
\author{Stephan Baier \and Arkaprava Bhandari} 

\address{Stephan Baier,
Ramakrishna Mission Vivekananda Educational and Research Institute, Department of Mathematics, G. T. Road, PO Belur Math, Howrah, West Bengal 711202, India}
\email{stephanbaier2017@gmail.com}
\address{Arkaprava Bhandari,
Ramakrishna Mission Vivekananda Educational and Research Institute, Department of Mathematics, G. T. Road, PO Belur Math, Howrah, West Bengal 711202, India}
\email{arkapravabhandari@gmail.com}


\begin{document} 
\begin{abstract}
Keating and Rudnick \cite{KR} derived asymptotic formulas for the variances of primes in arithmetic progressions and short intervals in the function field setting. Here we consider the hybrid problem of calculating the variance of primes in intersections of arithmetic progressions and short intervals. Keating and Rudnick used an involution to translate short intervals into arithmetic progressions. We follow their approach but apply this involution, in addition, to the arithmetic progressions. This creates dual arithmetic progressions in the case when the modulus $Q$ is a polynomial in $\mathbb{F}_q[T]$ such that $Q(0)\not=0$. The latter is a restriction which we keep throughout our paper. At the end, we discuss what is needed to relax this condition.
\end{abstract}
\maketitle
\tableofcontents

\section{Introduction to the problem}
In \cite{KR}, Keating and Rudnick investigated two important quantities from analytic number theory in the context of function fields: The variance of primes in arithmetic progressions to a fixed large modulus and the mean value of primes in short intervals, i.e. the quantities
$$
\int\limits_0^{N} \Bigg|\sum\limits_{x<n\le x+\Delta} \Lambda(n) -\Delta \Bigg|^2 dx
$$
and 
$$
\sum\limits_{\substack{a=1\\ (a,q)=1}}^q \Bigg|\sum\limits_{\substack{n\le N\\ n\equiv a\bmod{q}}} \Lambda(n)-\frac{N}{\varphi(q)} \Bigg|^2.
$$
The object of this paper is to establish a hybrid result on the variance of primes in  intersections of arithmetic progressions and short intervals in the function field setting. 

To formulate their results and our hybrid result, let us first set up some notations and recall fundamental relations, following \cite{KR}.

\section{Notations and fundamental relations}
\begin{itemize}
\item Let $\mathbb{F}_q$ be a field with $q=p^r$ elements ($p$ prime) and $\mathbb{F}_q[T]$ be the ring of polynomials in one variable over this field.
\item We define a norm of polynomials  $0\not=f\in \mathbb{F}_q[T]$ by 
$$
|f|=q^{\deg f}.
$$ 
\item We denote by $\mathcal{P}_n$ the set of polynomials of degree $n$, by $\mathcal{M}_n$ the subset of monic polynomials and by $\mathcal{P}_{\le n}$ the set of polynomials of degree less or equal $n$.
\item For $N\in \mathbb{F}_q[T]$ we define the von Mangoldt function as 
$$
\Lambda(N)=\begin{cases} \deg P & \mbox{ if } N=cP^k, \ P\in \mathcal{M}_n \mbox{ irreducible and } c\in \mathbb{F}_q^{\ast},\\ 
0 & \mbox{ otherwise.} \end{cases}
$$
We note the identity
\begin{equation} \label{primenumbertheorem}
\Psi(n):=\sum\limits_{N\in \mathcal{M}_n} \Lambda(N)=q^n,
\end{equation}
the prime number theorem for $\mathbb{F}_q[T]$.
\item If $C\in \mathcal{P}_n$ and $h<n$, we define
$$
I(C;h):=\{f : |f-C|\le q^h\}=C+\mathcal{P}_{\le h},
$$ 
which is the function field version of a short interval.
\item We define
$$
\nu(C;h):=\sum\limits_{\substack{f\in I(C;h)\\ f(0)\not=0}} \Lambda(f)
$$
and
$$
\Psi(n;Q,A):=\sum\limits_{\substack{N\in \mathcal{M}_n\\ N\equiv A\bmod{Q}}} \Lambda(N).
$$.
\item More generally, we define 
$$
\Psi(C,h;Q,A):=\sum\limits_{\substack{N\in I(C;h)\\N(0)\neq0\\N\equiv A~\text{mod}~ Q}}\Lambda(N),
$$
where we note that the summation is not restricted to monic polynomials. (Of course, if $C$ is monic and $h<\deg C$, then $N$ is forced to be monic.)
The restriction to $N(0)\neq 0$ excludes only polynomials of the form $N=cT^n$ since for all other polynomials with $N(0)=0$, the von Mangoldt function takes the value 0. This extra condition will help us in our calculations. 
\item The variance that we will be interested in is 
$$
V(n,h;Q):=\frac{1}{q^n}\sum\limits_{C\in \mathcal{M}_n} \sum\limits_{\substack{A \bmod{Q}\\ (A,Q)=1}} \left| \Psi(C,h;Q,A)-\frac{q^{h+1}}{\varphi(Q)}\right|^2.
$$
\item The Euler $\varphi$ function is defined as
$$
\varphi(Q):=\sharp\{A \bmod Q: (A,Q)=1\}
$$
and the M\"obius function $\mu$ is defined as 
$$
\mu(Q):=\begin{cases} (-1)^{\omega(Q)} & \mbox{ if } Q \mbox{ is square-free,} \\
0 & \mbox{ if } Q \mbox{ is not square-free,} \end{cases}
$$
where $\omega(Q)$ is the number of monic irreducible factors of $Q$.  
\item Assume $\deg Q\ge 2$ and let $\chi$ be a primitive Dirichlet character modulo $Q$. We call a character $\chi$ even if it is trivial on $\mathbb{F}_q^{\ast}$ and odd otherwise. Set
$$
\lambda_{\chi}:=\begin{cases} 1 & \mbox{ if } \chi \mbox{ is even,}\\ 0 & \mbox{ otherwise.} \end{cases}
$$
Set 
$$
d=\deg Q-1-\lambda_{\chi}.
$$
From the theory of Dirichlet $L$-function for function fields, it is known that there is a diagonal matrix
$$
\Theta_{\chi}=\mbox{diag}(e^{i\theta_1},...,e^{i\theta_d}),
$$
called unitarized Frobenius matrix,
having the property that the completed $L$-function belonging to $\chi$ satisfies 
$$
L^{\ast}(u,\chi)=\det(I-uq^{1/2}\Theta_{\chi})
$$
(see \cite[page 10]{KR}). Moreover, there is an explicit formula stating that
\begin{equation} \label{explicit}
\sum\limits_{N\in \mathcal{M}_n} \chi(n)\Lambda(n)=-q^{n/2}\text{tr }\Theta^n_{\chi}-\lambda_{\chi}.
\end{equation}
\item As in \cite[subsection 4.1.]{KR}, for a non-zero polynomial $X\in \mathbb{F}_q[T]$, we define its involution $X^{\ast}$ by $$X^{\ast}(T):=T^{\deg X}X\left(\frac{1}{T}\right). $$
\end{itemize}

\section{Results}
On the variance of primes in arithmetic progressions for function fields, Keating and Rudnick proved the following result.

\begin{theorem}[Theorem 2.2 in \cite{KR}] \label{KRresult} {\rm (i)} Fix $n\ge 1$. Given a finite field $\mathbb{F}_q$, let $Q(T)\in \mathbb{F}_q[T]$ be a polynomial such that $\deg Q>n$. Then
$$
\sum\limits_{\substack{A \bmod{Q}\\ (A,Q)=1}} \left| \Psi(n;Q,A)-\frac{q^n}{\varphi(Q)}\right|^2=nq^n-\frac{q^{2n}}{\varphi(Q)}+O\left(n^2q^{n/2}+(\deg Q)^2\right),
$$
the implied constant being absolute. 

{\rm (ii)} Fix $n\ge 2$. Given a sequence of finite fields $\mathbb{F}_q$ and square-free polynomials $Q(T)\in \mathbb{F}_q[T]$ of positive degree with $\deg Q\le n+1$, we have, as $q\rightarrow \infty$, 
$$
\sum\limits_{\substack{A \bmod{Q}\\ (A,Q)=1}} \left| \Psi(n;Q,A)-\frac{q^n}{\varphi(Q)}\right|^2\sim q^n(\deg Q-1).
$$ 
\end{theorem}

We note that there is an overlap of parts (i) and (ii): the situation when $\deg Q= n+1$. In this case, the result of part (i) implies that of part (ii) since $\varphi(Q)/q^{\deg Q}\rightarrow 1$ as $q\rightarrow \infty$. On the variance of primes in short intervals, Keating and Rudnick proved the following result.

\begin{theorem}[Theorem 2.1 in \cite{KR}] Fix $n\ge 4$ and $0\le h\le n-4$. Then, as $q\rightarrow \infty$, 
$$
\sum\limits_{C\in \mathcal{M}_n} \left|\nu(C;h)-q^{h+1}(1-q^{-n})\right|^2\sim q^{h+1}(n-h-2).
$$ 
\end{theorem}

It should be noted that Bank, Soroker and Rosenzweig \cite{BSR}, using advanced algebraic techniques, managed to prove the following asymptotic results on primes in {\it individual} short intervals and arithmetic progressions to large moduli. 

\begin{theorem}[Variants of Corollaries 2.4. and 2.5. in \cite{BSR}]  \label{BSRresult} {\rm (i)} For $n\in \mathbb{N}$ and $\varepsilon>0$ fixed, the asymptotic formula 
$$
\nu(C;h) \sim \sharp I(C;h), \quad \mbox{ as } q\rightarrow \infty
$$
holds uniformly for all $C\in \mathcal{M}_n$ and all $h\ge |C|^{\varepsilon}$. 

{\rm (ii)} For $n\in \mathbb{N}$, the asymptotic formula
$$
\Psi(n;Q;A)\sim \frac{\Psi(n)}{\varphi(q)}, \quad \mbox{ as } q\rightarrow \infty
$$
holds uniformly for all relative prime $Q,A\in \mathbb{F}_q[T]$ satisfying $|Q|\le q^{n(1-\delta)}$, where 
$$
\delta:=\begin{cases}
4/n & \mbox{ if } (A/Q)' \mbox{ is constant and } p=2,\\ 
3/n & \mbox{ otherwise.} \end{cases}
$$ 
\end{theorem}

Our main result is the following hybrid result on the variances considered by Keating and Rudnick. Part (iii) is conditional under a conjecture formulated in subsection \ref{Katz}.

\begin{theorem} \label{maintheorem} {\rm (i)}
Fix $n\ge 1$ and $0\le h\le n-1$. Given a finite field $\mathbb{F}_q$, let $Q(T)\in \mathbb{F}_q[T]$ be a polynomial such that $\deg Q>h$ and $Q(0)\not=0$. Then
\begin{equation*}
\begin{split}
V(n,h;Q)= nq^{h+1}- \frac{q^{2(h+1)}}{\varphi(Q)}+ O\left(\frac{n^2q^{h+1}}{q^{n/2}}+\frac{q^{h+1}(\deg Q)^2}{q^{n}}+\frac{q^{2(h+1)\deg Q}}{\varphi(Q)q^{n}}\right),
\end{split}
\end{equation*}
the implied constant being absolute. 

{\rm (ii)}  Fix $n\ge 1$ and $0\le h\le n-1$. Given a finite field $\mathbb{F}_q$, let $Q(T)\in \mathbb{F}_q[T]$ be a polynomial of positive degree such that $\deg Q\le n$ and $Q(0)\not=0$. Then
$$
V(n,h;Q)=\frac{q^{h+1}(q-1)}{\varphi(T^{n-h}Q^{\ast})} \sum\limits_{\substack{\chi \bmod T^{n-h}Q^{\ast}\\ \chi \ {\scriptsize \rm even\ primitive}}} |\text{\rm tr }\Theta^n_{\chi}|^2+O\left(q^{h}(n-h-1+\deg Q)^2\right), 
$$ 
the implied constant being absolute. In particular,
\begin{equation} \label{bound}
V(n,h;Q)\ll q^{h+1}(n-h-1+\deg Q)^2.
\end{equation}

{\rm (iii)} Fix $n\ge 5$ and $1\le h\le n-4$. Assume that Conjecture \ref{con} from subsection \ref{Katz} holds. Then given a sequence of finite fields $\mathbb{F}_q$ and square-free polynomials $Q(T)\in \mathbb{F}_q[T]$ with $3\le \deg Q\le h+2$ and $Q(0)\not=0$, we have, as $q\rightarrow \infty$,
$$
V(n,h;Q)\sim q^{h+1}(n-h-2+\deg Q).
$$ 
\end{theorem}

Again, there is an overlap of the parts (i) and (iii): the situations when $\deg Q=h+1$ or $\deg Q=h+2$. As $q\rightarrow \infty$, (i) implies that $V(n,h;Q)\sim q^{h+1}(n-1)$ if $\deg Q=h+1$ and $V(n,h;Q)\sim q^{h+1}n$ if $\deg Q=h+2$. Both agrees with $V(n,h;Q)\sim q^{h+1}(n-h-2+\deg Q)$ in these cases. This overlap provides a good test for our Conjecture \ref{con}.

Keating and Rudnick expressed the variance of primes in arithmetic progressions in terms of Dirichlet $L$-functions. To handle primes in short intervals, they used an involution to relate them to primes in arithmetic progressions, which is one of the key innovations in their paper. We follow their method. However, we need to take both into consideration: primes in short intervals {\it and} arithmetic progressions. This forces us to use the said involution on {\it both} sides. We therefore relate primes in the said arithmetic progressions to primes in {\it dual} arithmetic progressions. This will work if the modulus $Q$ satisfies $Q(0)\not=0$,
a condition which we will keep throughout this paper.  We will make some remarks on the case $Q(0)=0$ in the last section.  \\ \\ 
{\bf Acknowledgements.} The authors would like to thank the Ramakrishna Mission Vivekananda Educational and Research Institute for providing excellent working conditions. The second-named author would like to thank CSIR, Govt. of India for financial support in the form of a Junior Research Fellowship under file number 09/934(0015)/2019-EMR-I. We also thank the anonymous referee for his comments. 

\section{Calculation of a mean value}
We aim to evaluate the modified variance
\begin{equation} \label{var1}
\tilde{V}(n,h;Q):=\frac{1}{q^n}\sum_{C\in \mathcal{M}_n} \sum_{\substack{A~\text{mod}~ Q\\ (A,Q)=1}}\left| \Psi(C,h;Q,A)- \text{mean}(n,h;Q,A) \right|^2, 
\end{equation}
where $\text{mean}(n,h;Q,A)$ is a mean value of primes in arithmetic progressions, defined by 
\begin{equation} \label{meanvaluedef}
\text{mean}(n,h;Q,A)= \frac{1}{\varphi(Q)q^n}\sum_{C\in \mathcal{M}_n} \sum_{\substack{A~\text{mod}~ Q\\ (A,Q)=1}}\Psi(C,h;Q,A).
\end{equation}
First we calculate this mean value. We have 
\begin{equation} \label{meanvalue}
\begin{split}
    \text{mean}(n,h;Q,A)&= \frac{1}{\varphi(Q)q^{n-h-1}}\sum_{B\in \mathcal{M}_{n-h-1}} \sum_{\substack{A~\text{mod}~ Q\\ (A,Q)=1}}\sum_{\substack{N\in I(T^{h+1}B,h)\\N(0)\neq0\\N\equiv A~\text{mod}~ Q}}\Lambda(N)\\
    &= \frac{1}{\varphi(Q)q^{n-h-1}} \sum_{\substack{A~\text{mod}~ Q\\ (A,Q)=1}}\sum_{\substack{N\in \mathcal{M}_n\\ N\equiv A~\text{mod}~ Q\\ N(0)\not=0}}\Lambda(N)\\
    &= \frac{1}{\varphi(Q)q^{n-h-1}}\sum_{\substack{N\in \mathcal{M}_{n}\\(N,Q)= 1\\ N(0)\not=0}}\Lambda(N)\\
    &= \frac{1}{\varphi(Q)q^{n-h-1}}\Bigg(\sum_{N\in \mathcal{M}_{n}}\Lambda(N)-\sum_{\substack{P\in \mathcal{M}_n \text{ irreducible}\\ P|Q\\ \deg P|n}}\deg P-\Lambda(T^n)\Bigg)\\
    &=\frac{q^{h+1}}{\varphi(Q)}+O\left( \frac{q^{h+1}\deg Q}{\varphi(Q)q^{n}}\right),
\end{split}
\end{equation}
using the prime number theorem \eqref{primenumbertheorem}.
Therefore,
\begin{equation} \label{lastred}
\begin{split}
V(n,h;Q)- \tilde{V}(n,h;Q)
\ll & 
\frac{1}{q^n}\sum\limits_{C\in \mathcal{M}_n} \sum\limits_{\substack{A \bmod{Q}\\ (A,Q)=1}} \left( \Psi(C,h;Q,A)+\frac{q^{h+1}}{\varphi(Q)}\right) \cdot \frac{q^{h+1}\deg Q}{\varphi(Q)q^{n}}\\
\ll & \Bigg(\frac{q^{h+1}}{q^n}\sum\limits_{N\in \mathcal{M}_n} \Lambda(N)+ \frac{1}{q^n} \sum\limits_{C\in \mathcal{M}_n} \sum\limits_{\substack{A \bmod{Q}\\ (A,Q)=1}}\frac{q^{h+1}}{\varphi(Q)}\Bigg) \cdot \frac{q^{h+1}\deg Q}{\varphi(Q)q^{n}}
\\
\ll & \frac{q^{2(h+1)}\deg Q}{\varphi(Q)q^{n}},
\end{split}
\end{equation}
again using \eqref{primenumbertheorem}.

\section{Treatment of the range $\deg Q >h$ (Proof of Theorem \ref{maintheorem}(i))} 
Now we evaluate the variance in \eqref{var1} for the case when $\deg Q>h$.
In this case, the congruence $N\equiv A\bmod{Q}$ has at most one solution $N$ in $I(C;h)$ for every $C\in \mathcal{M}_n$.  Let $N(C,A)$ be this solution, if it exists. Then  
$$
\Psi(C,h;Q,A)=\begin{cases} \Lambda(N(C,A)), & \text{if } I(C;h)\cap (Q\mathbb{F}_q[T]+A)\setminus T\mathbb{F}_q[T]\not=\emptyset,\\ 0, & \text{otherwise}.
\end{cases}
$$
Hence, using \eqref{meanvalue}, the variance in question becomes 
\begin{equation*}
\begin{split}
\tilde{V}(n,h;Q) = & \frac{1}{q^n}\sum_{C\in \mathcal{M}_n} \sum_{\substack{A~\text{mod}~ Q\\ (A,Q)=1}}\left|\text{mean}(n,h;Q,A)-\begin{cases} \Lambda(N(C,A)), & \text{if } I(C;h)\cap (Q\mathbb{F}_q[T]+A)\setminus T\mathbb{F}_q[T]\not=\emptyset,\\ 0, & \text{otherwise}
\end{cases}\right|^2\\
	= & \frac{1}{q^n}\sum_{C\in \mathcal{M}_n}\Bigg(\sum_{\substack{N\in I(C;h)\\(N,Q)=1\\N(0)\neq0}}\Lambda(N)^2+\varphi(Q)\left(\frac{q^{h+1}}{\varphi(Q)}+ O\left( \frac{q^{h+1}\deg Q}{\varphi(Q)q^{n}}\right)\right)^2\Bigg)-\\
& -2\left(\frac{q^{h+1}}{\varphi(Q)}+O\left( \frac{q^{h+1}\deg Q}{\varphi(Q)q^{n}}\right)\right)\frac{1}{q^n}\sum_{C\in \mathcal{M}_n}\sum_{\substack{N\in I(C;h)\\(N,Q)=1\\N(0)\neq0}}\Lambda(N).
\end{split}
\end{equation*}
Now 
\begin{align*}
\sum_{C\in \mathcal{M}_n}\sum_{\substack{N\in I(C;h)\\(N,Q)=1\\N(0)\neq0}}\Lambda(N)&=q^{h+1}\sum_{\substack{N\in \mathcal{M}_n\\(N,Q)=1\\N(0)\neq0}}\Lambda(N)\\
&= q^{h+1}\Bigg(\sum_{N\in \mathcal{M}_n}\Lambda(N)-\sum_{\substack{P \text{ irreducible} \\P|Q \\ \deg P|n}}\deg P-\Lambda(T^n)\Bigg)\\
&= q^{h+1}(q^n+O(\deg Q))
\end{align*}
and
\begin{align*}
\sum_{C\in \mathcal{M}_n}\sum_{\substack{N\in I(C;h)\\(N,Q)=1\\N(0)\neq0}}\Lambda(N)^2 &= q^{h+1}\sum_{\substack{N\in \mathcal{M}_n\\(N,Q)=1\\N(0)\neq0}}\Lambda(N)^2\\
&= q^{h+1}\Bigg(\sum_{N\in \mathcal{M}_n}\Lambda(N)^2-\sum_{\substack{P \text{ irreducible}\\ P|Q \\ \deg P|n}}(\deg P)^2-\Lambda(T^n)^2 \Bigg)\\
&= q^{h+1}\left(nq^n+O\left(n^2q^{n/2}+(\deg Q)^2 \right)\right),
\end{align*}
again using \eqref{primenumbertheorem}. 
Thus the said variance becomes 
\begin{equation*}
\begin{split}
\tilde{V}(n,h;Q)= & nq^{h+1}- \frac{q^{2(h+1)}}{\varphi(Q)}+O\left(\frac{n^2q^{h+1}}{q^{n/2}}+\frac{q^{h+1}(\deg Q)^2}{q^{n}}+\frac{q^{2(h+1)}\deg Q}{\varphi(Q)q^n}\right).
\end{split}
\end{equation*}
This together with \eqref{lastred} proves Theorem \ref{maintheorem}(i). 

\section{Treatment of the range $\deg Q \le h+2$ (Proof of Theorem \ref{maintheorem}(ii))}
We recall that we defined the involution $X^{\ast}$ of a non-zero polynomial $X\in \mathbb{F}_q[T]$ by $$X^{\ast}(T):=T^{\deg X}X\left(\frac{1}{T}\right). $$
This involution reverses the order of the coefficients of $X$. 
It is straight-forward to check that it satisfies the following properties which will be used in the following. 
\begin{enumerate}
    \item $(XY)^{\ast}=X^{\ast}Y^{\ast}$ for all $X,Y\in \mathbb{F}_q[T]$.
    \item $(X+Y)^{\ast}=X^{\ast}+Y^{\ast}$ if deg $X=$deg $Y$.
    \item $(X^{\ast})^{\ast}=X$ iff $T\nmid X$.
    \item $\deg X=\deg X^{\ast}$ iff $T\nmid X$.
    \item $|X-BT^{h+1}|\leq q^h \Longleftrightarrow X^{\ast}\equiv B^{\ast}$ mod $T^{n-h}$ if $T\nmid X$, deg $X=n$ and deg $B=n-h-1$.
    \item $T\nmid B^{\ast}$ if $B\neq 0$.
    \item $\Lambda(X)=\Lambda(X^{\ast})$ if $T\nmid X$.
    \item $\varphi(X)=\varphi(X^{\ast})$ if $T\nmid X$.
\end{enumerate}
Now we assume $T\nmid XQ$ (or equivalently, $XQ(0)\not=0$) and $\Tilde{A}$ to be a representative of a given residue class $A \bmod Q$ such that $n=\deg X=\deg X^{\ast}=\deg \Tilde{A}$. (Such a representative $\tilde{A}$ surely exists if $\deg Q\le n$.) Then we have the following equivalences: \newpage
\begin{itemize}
    \item $X^{\ast}\equiv A$ mod $Q\Longleftrightarrow X^{\ast}-\Tilde{A}=VQ$ for some $V\in \mathbb{F}_q[T]$. \\{\bf Proof:} This is by the definition of congruences.
    \item $X^{\ast}-\Tilde{A}=VQ$ for some $V\in \mathbb{F}_q[T]\Longleftrightarrow (X^{\ast}-\Tilde{A})^{\ast}=V'Q^{\ast}$ for some $V'\in \mathbb{F}_q[T]$. \\ {\bf Proof of "$\Longrightarrow$":} Just take $V'=V^{\ast}$. \\ {\bf Proof of "$\Longleftarrow$":} Assume $T^w$ is the largest power of $T$ dividing $X^{\ast}-\tilde{A}$. Then 
$$
X^{\ast}-\tilde{A}=T^w(V'Q^{\ast})^{\ast}=T^w(V')^{\ast}Q
$$ 
since $T\nmid Q$. Now take $V=(V')^{\ast}T^w$. 
    \item $(X^{\ast}-\Tilde{A})^{\ast}=V'Q^{\ast}$ for some $V'\in \mathbb{F}_q[T]\Longleftrightarrow X-\Tilde{A}^{\ast}=V'Q^{\ast}$ for some $V'\in \mathbb{F}_q[T]$. \\ {\bf Proof:} Since $\deg X^{\ast}=\deg \tilde{A}$, we have $(X^{\ast}-\tilde{A})^{\ast}=(X^{\ast})^{\ast}-\tilde{A}^{\ast}=X-\tilde{A}^{\ast}$, where the last equation holds because of $T\nmid X$.
    \item $X-\Tilde{A}^{\ast}=V'Q^{\ast}$ for some $V'\in \mathbb{F}_q[T]\Longleftrightarrow X \equiv \Tilde{A}^{\ast}$ mod $Q^{\ast}$. \\ {\bf Proof:} This is by the definition of congruences.  
\end{itemize}
(We note that the condition $T\nmid Q$ is crucial in the implication "$\Longleftarrow$" of the second equivalence.)  
So altogether $X^{\ast}\equiv A \bmod Q \Longleftrightarrow X\equiv \Tilde{A}^{\ast} \bmod Q^{\ast}$ under the above conditions. 

Throughout the following, we keep the conditions that $T\nmid Q$, $(A,Q)=1$, $\deg Q\le n$ and $B\in \mathcal{P}_{n-h-1}\setminus \{0\}$. Then, as a priliminary step, we deduce that
\begin{equation*}
\begin{split}
\Psi(T^{h+1}B,h;Q,A)= & \sum_{\substack{N\in I(T^{h+1}B;h)\\N(0)\neq0\\N \equiv A~\text{mod}~ Q}}\Lambda(N)=\sum_{\substack{\deg N^{\ast}=n\\ N^{\ast}\equiv B^{\ast}\text{ mod }T^{n-h}\\ N\equiv A~\text{mod}~ Q}}\Lambda(N)\\ = &\sum_{\substack{\deg N=n\\ N\equiv B^{\ast}\text{ mod }T^{n-h}\\ N^{\ast}\equiv A~\text{mod}~ Q}}\Lambda(N)=\sum_{\substack{\deg N=n\\ N\equiv B^{\ast}\text{ mod }T^{n-h}\\ N\equiv \Tilde{A}^{\ast}~\text{mod}~ Q^{\ast}}}\Lambda(N),
\end{split}
\end{equation*}
where we recall that $N$ is not restricted to monic polynomials. (Of course, if $B$ is monic, then $N$ is forced to be monic.)
In this way, we have rewritten $\Psi(T^{h+1}B,h;Q,A)$ as sum of the von Mangoldt function over arithmetic progressions. This can be treated via character sums over primes and subsequently using Dirichlet $L$-functions by standard techniques. We further observe that $(\Tilde{A}^{\ast},Q^{\ast})=1$ if $(A,Q)=1$ and also $(B^{\ast},T^{n-h})=1$, the latter being trivial. To see $(\tilde{A}^{\ast},Q^{\ast})=1$, assume that $\tilde{A}^{\ast}=ED$ and $Q^{\ast}=FD$ for some $D,E,F\in \mathbb{F}_q[T]$. Now by multiplicativity, $(\tilde{A}^{\ast})^{\ast}=E^{\ast}D^{\ast}$ and $Q=(Q^{\ast})^{\ast}=F^{\ast}D^{\ast}$. But $(\tilde{A}^{\ast})^{\ast}T^w=\tilde{A}$ for some $w\in \mathbb{N}_0$. So $D^{\ast}|\tilde{A}$ and $D^{\ast}|Q$. This forces $D^{\ast}\in \mathbb{F}_q^{\ast}$ since $(\tilde{A},Q)=1$, which implies $D=cT^v$ for some $v\in \mathbb{N}_0$ and $c\in \mathbb{F}_q^{\ast}$. However, $T\nmid Q^{\ast}$, so $D\in \mathbb{F}_q^{\ast}$, proving $(\tilde{A}^{\ast},Q^{\ast})=1$. Moreover, since $\varphi(Q)=\varphi(Q^{\ast})$ and $\tilde{A}_1^{\ast}- \tilde{A}_2^{\ast} \equiv 0 \bmod{Q^{\ast}}$ implies $A_1-A_2\equiv \tilde{A}_1-\tilde{A}_2\equiv 0\bmod Q$, it follows that if $A$ runs over all residue classes modulo $Q$ with $(A,Q)=1$, then $\tilde{A}^{\ast}$ runs over all residue classes $A'$ modulo $Q^{\ast}$ with $(A',Q^{\ast})=1$. 

Using the above considerations and the orthogonality relations for Dirichlet characters $\chi \bmod{T^{n-h}}$ and $\chi' \bmod{Q^{\ast}}$, we obtain
\begin{equation} \label{shortintervalsum}
\Psi(T^{h+1}B,h;Q,A)=\frac{1}{\varphi(T^{n-h})\varphi(Q^{\ast})}\sum_{\substack{\chi \text{ mod }T^{n-h}\\\chi'\text{ mod }Q^{\ast}}}\overline{\chi}(B^{\ast})\overline{\chi'}(\Tilde{A}^{\ast})\sum_{\deg N=n}\chi\chi'(N)\Lambda(N). 
\end{equation}
Let $\Tilde{\chi}=\chi\chi'$ and $\Tilde{Q}=T^{n-h}Q^{\ast}$. We note that 
$$
\varphi(\Tilde{Q})=\varphi(T^{n-h})\varphi(Q^{\ast})=\varphi(T^{n-h})\varphi(Q)=(q-1)q^{n-h-1}\varphi(Q).
$$ 
Now using \eqref{meanvalue}, $\Lambda(N)=\Lambda(N^{\ast})$ and $(N,Q)=1=(N^{\ast},Q^{\ast})$, we observe that the contribution of the trivial character $\Tilde{\chi}=\Tilde{\chi}_0$ to the right-hand side equals exactly the mean value $\text{mean}(n,h;Q,A)$. It follows that
\begin{equation*} 
\begin{split}
  \tilde{V}(n,h;Q) = &\frac{1}{q^n}\sum_{C\in \mathcal{M}_n} \sum_{\substack{A~\text{mod}~ Q\\ (A,Q)=1}}\left|\Psi(C,h;Q,A)- \text{mean}(n,h;Q,A) \right|^2\\
    =~& \frac{1}{q^{n-h-1}}\sum_{B\in \mathcal{M}_{n-h-1}} \sum_{\substack{A~\text{mod}~ Q\\ (A,Q)=1}}\left|\Psi(T^{h+1}B,h;Q,A)- \text{mean}(n,h;Q,A) \right|^2\\
=~& \frac{1}{(q-1)q^{n-h-1}}\sum_{B\in \mathcal{P}_{n-h-1}} \sum_{\substack{A~\text{mod}~ Q\\ (A,Q)=1}}\left|\Psi(T^{h+1}B,h;Q,A)- \text{mean}(n,h;Q,A) \right|^2\\
    =~& \frac{1}{(q-1)q^{n-h-1}}\sum_{B\in \mathcal{P}_{n-h-1}} \sum_{\substack{A~\text{mod}~ Q\\ (A,Q)=1}}\Big|\frac{1}{\varphi(T^{n-h})\varphi(Q^{\ast})}\sum_{\substack{\chi \text{ mod }T^{n-h}\\\chi'\text{ mod }Q^{\ast}}}\overline{\chi}(B^{\ast})\overline{\chi'}(\Tilde{A}^{\ast})\sum_{\deg N=n}\chi\chi'(N)\Lambda(N)\\
& - \text{mean}(n,h;Q,A) \Big|^2\\
    =~& \frac{1}{(q-1)q^{n-h-1}}\sum_{\substack{B^{\ast}\text{ mod }T^{n-h}\\ (B^{\ast},T)=1}} \sum_{\substack{\tilde{A}^{\ast}~\text{mod}~ Q^{\ast}\\ (A,Q)=1}}\Big|\frac{1}{\varphi(T^{n-h})\varphi(Q^{\ast})}\sum_{\substack{\chi \text{ mod }T^{n-h}\\\chi'\text{ mod }Q^{\ast}\\ \chi\chi'\neq \Tilde{\chi}_0}}\overline{\chi}(B^{\ast})\overline{\chi'}(\Tilde{A}^{\ast})\sum_{\deg N=n}\chi\chi'(N)\Lambda(N)\Big|^2\\
 =~& \frac{1}{(q-1)q^{n-h-1}\varphi(\Tilde{Q})}\sum_{\substack{\Tilde{\chi} \text{ mod }\Tilde{Q}\\ \Tilde{\chi}\neq\Tilde{\chi}_0}}\left|\Psi(n,\Tilde{\chi}) \right|^2\\
=~& \frac{1}{(q-1)q^{n-h-1}\varphi(\Tilde{Q})}\sum_{\substack{\Tilde{\chi} \text{ mod }\Tilde{Q} \text{ even,}\\ \Tilde{\chi}\neq\Tilde{\chi}_0}}\left|\Psi(n,\Tilde{\chi}) \right|^2,
\end{split}
\end{equation*}
where 
$$
\Psi(n,\Tilde{\chi})=\sum_{\deg N=n}\Tilde{\chi}(N)\Lambda(N).
$$ 
Below we include some explanations of several steps. For the third equation, we observe that $\mathcal{P}_{n-h-1}$ can be written as a union $\bigcup_{c\in \mathbb{F}_q^{\ast}} c\mathcal{M}_{n-h-1}$. If $B\in \mathcal{M}_{n-h-1}$ and $c\in \mathbb{F}_q^{\ast}$, then $\Psi(T^{h+1}B,h;Q,A)=\Psi(T^{h+1}cB,h;Q,cA)$, and if $A$ runs over all residue classes modulo $Q$ with $(A,Q)=1$, then $cA$ does as well. 
For the second-last equation, we again use the orthogonality relations for Dirichlet characters after squaring out. 
For the last equation, we note that if $\Tilde{\chi}$ is odd, then $\psi(n,\Tilde{\chi})$ vanishes because in this case $\sum_{c\in \mathbb{F}_q^{\ast}} \Tilde{\chi}(cN)$ equals 0 for every $N$. 

Now, following Keating and Rudnick \cite{KR}, we break the sum over even characters $\tilde{\chi}$ mod $\tilde{Q}$ into sums over primitive even and non-primitive even characters. Using \cite[(3.22)]{KR}, the number of primitive even characters modulo $\tilde{Q}$ equals
$$
\varphi_{\text{prim}}^{\text{ev}}(\tilde{Q})=\frac{1}{q-1}\sum\limits_{D|\tilde{Q}} \mu(D)\varphi\left(\frac{\tilde{Q}}{D}\right)=\frac{\varphi(\varphi(\tilde{Q}))}{q-1}. 
$$
Since the number of even characters modulo $\tilde{Q}$ equals 
$$
\varphi^{\text{ev}}(\tilde{Q})=\frac{\varphi(\tilde{Q})}{q-1},
$$
the number of non-primitive even characters modulo $\tilde{Q}$ therefore equals
$$
\varphi_{\text{non-prim}}^{\text{ev}}(\tilde{Q})=\frac{\varphi(\tilde{Q})-\varphi(\varphi(\tilde{Q}))}{q-1}.
$$ 
As $q\rightarrow\infty$, it follows that (cf. \cite[subsection 3.3.]{KR})
$$
\varphi_{\text{non-prim}}^{\text{ev}}(\tilde{Q})\ll
 \frac{\varphi(\tilde{Q})}{q(q-1)}=\frac{\varphi^{\text{ev}}(\tilde{Q})}{q}
$$
and 
\begin{equation} \label{primeven}
\varphi_{\text{prim}}^{\text{ev}}(\tilde{Q})=\varphi^{\text{ev}}(\tilde{Q})\left(1+O\left(\frac{1}{q}\right)\right)=\frac{\varphi(\tilde{Q})}{q-1}\left(1+O\left(\frac{1}{q}\right)\right). 
\end{equation}
Now for non-primitive and non-trivial characters $\tilde{\chi}$ mod $\tilde{Q}$, we use the bound 
$$
\left|\sum\limits_{N\in \mathcal{M}_n} \tilde{\chi}(N)\Lambda(N)\right|\le q^{n/2}(\deg \tilde{Q}-1)=q^{n/2}(n-h-1+\deg Q)
$$ 
following from the Riemann hypothesis (see \cite[(3.33)]{KR}).  Since for even characters $\tilde{\chi}$, 
\begin{equation} \label{inflation}
\Psi(n,\tilde{\chi})=(q-1)\sum\limits_{N\in \mathcal{M}_n} \tilde{\chi}(N)\Lambda(N), 
\end{equation}
we deduce that
\begin{equation*}
\tilde{V}(n,h;Q) =  \frac{1}{(q-1)q^{n-h-1}\varphi(\Tilde{Q})} \sum_{\substack{\Tilde{\chi} \text{ mod }\Tilde{Q}\\\Tilde{\chi} \text{ even\ primitive }}} \left|\Psi(n,\Tilde{\chi}) \right|^2 
+O\left(q^{h}(n-h-1+\deg Q)^2\right),
\end{equation*}
the $O$-term accounting for the non-primitive even characters. 
Taking \eqref{inflation} into account,
the explicit formula \eqref{explicit} now implies
\begin{equation} \label{Theta}
\tilde{V}(n,h;Q)
= \frac{q^{h+1}(q-1)}{\varphi(\tilde{Q})} \sum\limits_{\substack{\tilde{\chi} \bmod \tilde{Q}\\ \tilde{\chi} \ {\scriptsize \rm even\ primitive}}} |\text{tr }\Theta^n_{\tilde{\chi}}|^2+O\left(q^{h}(n-h-1+\deg Q)^2\right),
\end{equation}
which together with \eqref{lastred} completes the proof of Theorem \ref{maintheorem}(ii), where the bound \eqref{bound} follows again using the Riemann hypothesis and \eqref{primeven}.

\section{Equidistribution (Proof of Theorem \ref{maintheorem}(iii))} \label{Katz}
Keating and Rudnick \cite{KR} invoked the following results of Katz from \cite{Ka1} and \cite{Ka2} (\cite[Theorems 4.2 and 5.2]{KR}) to obtain asymptotics as $q\rightarrow \infty$.

\begin{theorem} \label{Katztheorem}
{\rm (i)} Fix $l\ge 4$. The conjugacy classes of $\Theta_{\chi}$ for the family of even primitive characters $\chi$ modulo $T^{l}$ become equidistributed in the projective unitary group $PU(l-2)$ of size $l-2$, as $q\rightarrow \infty$. 

{\rm (ii)} Fix $m\ge 2$. Suppose we are given a sequence of fields $\mathbb{F}_q$ and square-free polynomials $Q(T)\in \mathbb{F}_q[T]$ of degree $m$. As $q\rightarrow \infty$, the conjugacy classes of $\Theta_{\chi}$ with $\chi$ running over all primitive odd characters modulo $Q$, become  equidistributed in the unitary group $U(m-1)$. 
\end{theorem}

As consequences, they deduced
$$
\lim\limits_{q\rightarrow \infty} \frac{1}{\phi^{\text{ev}}_{\text{prim}}(T^{n-h})} \sum\limits_{\substack{\chi \bmod{T^{n-h}}\\ \chi \text{ primitive even }}} \left|\text{tr } \Theta_{\chi}^n\right|^2 = \int\limits_{\text{PU}(n-h-2)} \left|\text{tr } U^n\right|^2dU = \min\{n,n-h-2\}=n-h-2
$$
if $0\le h\le n-4$ and 
$$
\lim\limits_{q\rightarrow \infty} \frac{1}{\phi^{\text{odd}}_{\text{prim}}(Q)} \sum\limits_{\substack{\chi \bmod{Q}\\ \chi \text{ primitive odd }}} \left|\text{tr } \Theta_{\chi}^n\right|^2 = \int\limits_{\text{U}(\deg Q -1) } \left|\text{tr } U^n\right|^2dU = \min\{n,\deg Q-1\}=\deg Q -1
$$  
if $2\le \deg Q\le n+1$, 
where $\phi^{\text{odd}}_{\text{prim}}(Q)$ is the number of primitive odd characters modulo $Q$ (for $Q$ square-free and of degree exceeding 1). In \cite{Ka1}, Katz also investigated the distribution of $\Theta_{\chi}$ for primitive even characters to square-free moduli $Q$ of degree exceeding $2$ and gave some partial results and a conjecture in this direction. Here we cautiously pose the following conjecture on the second moment of traces of Frobenii for even characters to moduli of the form $T^lQ^m$ which seems to consist with these investigations.  

\begin{conjecture} \label{con} Suppose $l\ge 4$ and $m\ge 3$. Suppose we are given a sequence of fields $\mathbb{F}_q$ and square-free polynomials $Q(T)\in \mathbb{F}_q[T]$ of degree $m$ and coprime to $T$. Then, as $q\rightarrow \infty$,  we have 
$$
\lim\limits_{q\rightarrow \infty} \frac{1}{\phi^{\text{\rm ev}}_{\text{\rm prim}}(T^lQ)} \sum\limits_{\substack{\chi \bmod{T^lQ}\\ \chi \ \text{\rm primitive even }}} \left|\text{\rm tr } \Theta_{\chi}^n\right|^2 = \min\{n,l+m-2\}.
$$
\end{conjecture}

This would imply
\begin{equation*}
\lim\limits_{q\rightarrow \infty} \frac{1}{\phi^{\text{ev}}_{\text{prim}}(\tilde{Q})} \sum\limits_{\substack{\tilde{\chi} \bmod{\tilde{Q}}\\ \tilde{\chi} \text{ primitive even }}} \left|\text{tr } \Theta_{\tilde{\chi}}^n\right|^2 =\min\{n,n-h-2+\deg Q\}=n-h-2+\deg Q
\end{equation*}
if $1\le h\le n-4$ and $3\le \deg Q\le h+2$. Combining this with \eqref{primeven} and \eqref{Theta}, we would get
$$
\tilde{V}(n,h;Q)\sim q^{h+1}(n-h-2+\deg Q),
$$
which together with \eqref{lastred} establishes the result of Theorem \ref{maintheorem}(iii).  

\section{Remarks on the case $Q(0)=0$}
If $Q(0)=0$, then our above method of converting the congruence 
$N\equiv A\bmod{Q}$ into a dual congruence $N^{\ast} \equiv \tilde{A}^{\ast} \bmod{Q^{\ast}}$ fails, as observed above. An extreme case is when $Q=T^m$ for some $m$. In this case, the involution will translate the congruence $N\equiv A\bmod{Q}$ back into a short interval condition on $N^{\ast}$, which is exactly what we aimed to get rid of. So in the general case when $Q(0)=0$ is allowed, one will ultimately need to deal with sums of the form
$$
\sum\limits_{\substack{N \ \mbox{\scriptsize \rm monic}\\ N(0)\not=0\\ \deg N = n =\deg N^{\ast}}} \chi(N)\chi^{\ast}(N^{\ast})\Lambda(N),
$$
where $\chi$ and $\chi^{\ast}$ are characters to different moduli. In fact, it will suffice to consider pairs of characters $(\chi,\chi^{\ast})$, where $\chi^{\ast}$ is a character modulo a power $T^m$. The above character sums relate to generalized $L$-functions of the form
$$
L(s,\chi,\chi^{\ast})=\sum\limits_{\substack{N\in \mathbb{F}_q[T]\\ N(0)\not=0\\ N \mbox{\scriptsize \rm \ monic}}} \chi(N)\chi^{\ast}(N^{\ast})\mathcal{N}(N)^{-s}.
$$
Note that the function 
$$
F(N):=\chi(N)\chi^{\ast}(N^{\ast})
$$
is completely multiplicative on the set of monic polynomials $N$ coprime to $T$. Therefore, $L(s,\chi,\chi^{\ast})$ can be written as an Euler product
$$
L(s,\chi,\chi^{\ast})=\prod\limits_{\substack{P \mbox{\scriptsize\rm \ monic,  irreducible}\\ P\not=T}} \left(1-\chi(P)\chi^{\ast}(P^{\ast})\mathcal{N}(P)^{-s}\right)^{-1}.
$$ 
It needs to be investigated which analytic properties these generalized $L$-functions have. (Are they rational functions? Does a Riemann Hypothesis hold for them?)

\end{document}